\documentclass[a4paper,12pt]{article}
\usepackage{amsfonts, amsmath, amssymb, amscd, latexsym, graphicx}
\newtheorem{thm}{Theorem}[section]
\newtheorem{prop}[thm]{Proposition}
\newtheorem{lem}[thm]{Lemma}
\newtheorem{cor}[thm]{Corollary}

\newtheorem{defn}[thm]{Definition}

\title{\Large Abstract commensurators of solvable Baumslag -- Solitar groups}

\author{O.\ Bogopolski \\ \small{Institute of Mathematics of} \\ \small{Siberian Branch of Russian Academy
of Sciences,} \\ {\small Novosibirsk, Russia}\\ {\small and D\"{u}sseldorf University, Germany} \\ \small{e-mail:
Oleg$\_$Bogopolski@yahoo.com}}

\large

\begin{document}

\maketitle

\begin{abstract} \noindent We prove that for any $n\geqslant 2$, the abstract commensurator group of the Baumslag -- Solitar group ${\text{\rm BS}}(1,n)$ is isomorphic to the subgroup $\{ \begin{pmatrix}
1 & q\\
0 & p
\end{pmatrix}
\,|\, q\in \mathbb{Q},\, p\in \mathbb{Q}^{\ast}
\}$
of ${\text{\rm GL}}_2(\mathbb{Q})$.

We also prove that for any finitely generated group $G$ with the unique root property
the natural homomorphisms ${\text{\rm Aut}}(G)\rightarrow {\text{\rm Comm}}(G)\rightarrow
{\text{\rm QI}}(G)$ are embeddings.

\end{abstract}

\section{Introduction}

For a group $G$, we denote by ${\text{\rm Aut}}(G)$ its automorphism group, by ${\text{\rm Comm}}(G)$ its abstract commensurator group, and by ${\text{\rm QI}}(G)$ its quasi-isometry group;
see Definitions~\ref{2.1} and~\ref{5.1}.
For a finitely generated $G$, there are natural homomorphisms
$${\text{\rm Aut}}(G)\rightarrow {\text{\rm Comm}}(G)\rightarrow
{\text{\rm QI}}(G),$$
which became embeddings if $G$ has the unique root property,
i.e. if
$$\forall x,y\in G\,\, \forall n\in \mathbb{N}\,\, (x^n=y^n \Rightarrow x=y);$$
see Sections~2 and~5.

We are interested in computing of abstract commensurator groups of
(solvable) Baumslag -- Solitar groups.
 The Baumslag -- Solitar groups ${\text{\rm BS}}(m,n)$,  $1\leqslant m\leqslant n,$ are given by the presentation $\langle a,b\,|\, a^{-1}b^ma=b^n\rangle$. These groups have served as a proving ground for many new ideas in
combinatorial and geometric group theory (see, for instance,~\cite{BS,FM,W}).
The only solvable groups in this class are groups ${\text{\rm BS}}(1,n)$; the groups ${\text{\rm BS}}(m,n)$ with
$1<m\leqslant n$ contain a free nonabelian group.

The automorphism groups of ${\text{\rm BS}}(m,n)$ were described by Collins in~\cite{C}.
It follows that the automorphism groups of ${\text{\rm BS}}(1,n)$ and ${\text{\rm BS}}(1,k)$ with $n,k\geqslant 1$ are isomorphic if and only if $n$ and $k$ have the same sets of prime divisors.

In~\cite{FM}, Farb and Mosher proved for $n\geqslant 2$ that ${\text{\rm QI}}({\text{\rm BS}}(1,n))\cong {\text{\rm Bilip}}(\mathbb{R})\times  {\text{\rm Bilip}}(\mathbb{Q}_n)$,
where $\mathbb{Q}_n$ is the metric space of $n$-adic rationals with the usual metric and ${\text{\rm Bilip}}(Y)$ denotes the group of bilipschitz homeomorphisms of a
metric space $Y$.

\noindent Moreover, they proved that ${\text{\rm BS}}(1,n)$
and ${\text{\rm BS}}(1,k)$ with $n,k\geqslant 1$ are quasi-isometric if and only if these groups are commensurable, that happens if and only if $n$ and $k$ have common powers.
In~\cite{W}, Whyte proved that groups ${\text{\rm BS}}(m,n)$ with $1 <
m < n$ are quasi-isometric.

In this paper we compute the abstract commensurator groups of ${\text{\rm BS}}(1,n)$.
It turns out that the abstract commensurator groups of all groups ${\text{\rm BS}}(1,n)$, $n\geqslant 2$,
are isomorphic.


\medskip

{\bf Main Theorem.} (Theorem~\ref{4.5}) {\it For every $n\geqslant 2$, ${\text{\rm Comm}}({\text{\rm BS}}(1,n))$ is isomorphic
to the subgroup $\{ \begin{pmatrix}
1 & q\\
0 & p
\end{pmatrix}
\,|\, q\in \mathbb{Q},\, p\in \mathbb{Q}^{\ast}
\}$
of ${\text{\rm GL}}_2(\mathbb{Q})$.}

\medskip

Note that ${\text{\rm BS}}(1,1)\cong \mathbb{Z}^2$, and it is well known that ${\text{\rm Comm}}(\mathbb{Z}^m)\cong
{\text{\rm GL}}_m(\mathbb{Q})$ for $m\geqslant 1$.

\section{General facts on commensurators}

\begin{defn}\label{2.1}
{\rm Let $G$ be a group. Consider the set $\Omega (G)$ of all isomorphisms
between subgroups of finite index of $G$. Two such isomorphisms
$\varphi_1:H_1\to H_1'$ and $\varphi_2:H_2\to H_2'$ are called
\emph{equivalent}, written $\varphi_1\sim\varphi_2$, if there exists a
subgroup $H$ of finite index in $G$ such that both $\varphi_1$ and
$\varphi_2$ are defined on $H$ and $\varphi_1|_H=\varphi_2|_H$.

For any two isomorphisms $\alpha:G_1\to G_1'$ and $\beta:G_2\to G_2'$
in $\Omega(G)$, we define their product
$\alpha\beta:\alpha^{-1}(G_1'\cap G_2)\to \beta(G_1'\cap G_2)$ in
$\Omega (G)$. The factor-set $\Omega(G)/{\sim}$ inherits the
multiplication $[\alpha][\beta]=[\alpha\beta]$ and is a group, called
the \emph{abstract commensurator group} of $G$ and denoted ${\text{\rm Comm}}(G)$.}
\end{defn}

\begin{defn} {\rm A group $G$ has the {\it unique root property} if for any $x,y\in G$ and any positive integer
$n$, the equality $x^n=y^n$ implies $x=y$.}
\end{defn}

For closeness, we reproduce here short proofs of the following two statements from~\cite{BB}.

\begin{prop}\label{2.3}
  Let $G$ be a group with the unique root property. Then ${\text{\rm Aut}}(G)$
  naturally embeds in ${\text{\rm Comm}}(G)$.
\end{prop}

{\it Proof.}
There is a natural homomorphism ${\text{\rm Aut}}(G)\rightarrow
{\text{\rm Comm}}(G)$. Suppose that some $\alpha\in {\text{\rm Aut}}(G)$ lies in its kernel.
Then $\alpha|_H={\text{\rm id}}$ for some subgroup $H$ of finite index in $G$.
If $m$ is this index, then $g^{m!}\in H$ for every $g\in G$.
Then $\alpha(g^{m!})=g^{m!}$. Extracting roots, we get $\alpha(g)=g$, that is $\alpha={\text{\rm id}}$. \hfill $\Box$

\begin{lem}\label{2.4}
Let $G$ be a group with the unique root property. Let $\varphi_1: H_1\rightarrow H_1'$ and
$\varphi_2: H_2\rightarrow H_2'$ be two isomorphisms between subgroups of finite index in $G$. Suppose
that $[\varphi_1]=[\varphi_2]$ in ${\text{\rm Comm}}(G)$. Then $\varphi_1|_{H_1\cap H_2}=\varphi_2|_{H_1\cap H_2}$.
\hfill $\Box$
\end{lem}

{\it Proof.}
The equality $[\varphi_1]=[\varphi_2]$ means that there exists a
  subgroup $H$ of finite index in $G$ such that both $\varphi_1$ and
  $\varphi_2$ are defined on $H$ and
  $\varphi_1|_H=\varphi_2|_H$. Clearly $H\leqslant H_1\cap
  H_2$. Denote $m=|(H_1\cap H_2):H|$.  Let $h$ be an arbitrary element
  of $H_1\cap H_2$. Then $h^{m!}\in H$ and so
  $\varphi_1(h^{m!})=\varphi_2(h^{m!})$. Since $G$ is a group with the
  unique root property, we get $\varphi_1(h)=\varphi_2(h)$. \hfill $\Box$

\begin{lem}\label{2.5} The group ${\text{\rm BS}}(1,n)$ has the unique root property.
In particular, ${\text{\rm Aut}}({\text{\rm BS}}(1,n))$ naturally embeds in ${\text{\rm Comm}}({\text{\rm BS}}(1,n))$.
\end{lem}

{\it Proof.}
The first claim follows by using matrix calculations in view of Lemma~\ref{4.1}. The second claim follows from Proposition~\ref{2.3}. \hfill $\Box$

\section{A structure of finite index subgroups of ${\text{\rm BS}}(1,n)$}

Let ${\text{\rm BS}}(1,n)=\langle a,b\,|\, a^{-1}ba=b^n\rangle$, where $n\geqslant 2$.
Denote $b_j=a^{-j}ba^j$,  $j\in \mathbb{Z}$. Then $$b_j^n=b_{j+1},\hspace*{4mm} a^{-1}b_ja=b_{j+1},\hspace*{4mm} b_ib_j=b_jb_i\hspace*{4mm} (i,j\in \mathbb{Z}).$$
Consider the homomorphism
$$\begin{array}{lr}
\psi: & {\text{\rm BS}}(1,n)\rightarrow \mathbb{Z}\\
& a\mapsto 1\\
& b\mapsto 0.\hspace*{-1mm}
\end{array}$$

\begin{lem}\label{3.1} {\rm 1)} We have ${\text{\rm BS}}(1,n)=U\rtimes V$, where $U={\text{\rm ker}}\,\psi=\langle b_j\,|\, j\in \mathbb{Z}\rangle$, $V=\langle a\rangle $, and $V$ acts on $U$ by the rule $a^{-1}b_ja=b_{j+1}$.

{\rm 2)} The subgroup $U$ has the presentation
$\langle b_j\,|\, b_j^n=b_{j+1},\hspace*{2mm}j\in \mathbb{Z}\rangle$ and so it can be identified with $\mathbb{Z}[\frac{1}{n}]$.

{\rm 3)} ${\text{\rm BS}}(1,n)\cong \mathbb{Z}[\frac{1}{n}]\rtimes \mathbb{Z}$, where $\mathbb{Z}$ acts on $\mathbb{Z}[\frac{1}{n}]$ by multiplication by $n$.
\end{lem}

{\it Proof.} The first claim is obvious, the second follows by applying the Reidemeister -- Schreier method,
and the third claim follows from the first two. \hfill $\Box$

\begin{lem}\label{3.2} Every subgroup $H$ of finite index in ${\text{\rm BS}}(1,n)$ can be written as
$H=\langle a^ku,w\rangle$ for some $k> 0$, $u,w\in U$ and $w\neq 1$.
\end{lem}

{\it Proof.} The subgroup $H$ is finitely generated.
Since the image of $H$ under the epimorphism $\psi:
{\text{\rm BS}}(1,n)\rightarrow \mathbb{Z}$
is generated by some $k>0$,
we can write $H=\langle a^ku,u_1,\dots ,u_s\rangle$ for some $u,u_1,\dots ,u_s\in U={\text{\rm ker}}\,\psi$. Observe that every finitely generated subgroup of $U\cong \mathbb{Z}[\frac{1}{n}]$ is cyclic. So, $H=\langle a^ku,w\rangle$ for some
$w\in U$. Clearly, $w\neq 1$, otherwise ${\text{\rm BS}}(1,n)$ were virtually cyclic, that is impossible. \hfill $\Box$

\begin{lem}\label{3.3} Let $H=\langle a^kb_q^r,b_p^s\rangle$ with $k>0$. Then $H=\langle a^kb_q^r,b_i^s\rangle$ for every $i\in\mathbb{Z}$.
\end{lem}

{\it Proof.} Since $(a^kb_q^r)^{-t}\cdot b_p^s\cdot (a^kb_q^r)^t=b_{p+tk}^s$ for every integer $t$, we have
$$H=\langle a^kb_q^r,b_{p+tk}^s\rangle=\langle a^kb_q^r,b_{p+(t+1)k}^s\rangle.$$
Given $i\in \mathbb{Z}$, we choose $t$ such that $p+tk\leqslant i<p+(t+1)k$. Then
$H=\langle a^kb_q^r,b_i^s\rangle$, since $b_i$ is a power of $b_{p+tk}$ and $b_{p+(t+1)k}$ is a power of $b_i$. \hfill $\Box$

\begin{prop}\label{3.4} Every subgroup $H$ of finite index in ${\text{\rm BS}}(1,n)$ can be written as
$H=\langle a^kb^l,b^m\rangle$ for some integer $k,l,m$, where $k,m>0$ and $(m,n)=1$.
The index of this subgroup is $km$.
\end{prop}

{\it Proof.} By Lemma~\ref{3.2}, $H=\langle a^kb_q^r,b_p^s\rangle$ for some $k,s> 0$ and $r,q,p\in \mathbb{Z}$. Set $m=s/(n,s)$. Clearly, $(m,n)=1$. We claim that $H=\langle a^kb_q^r,b_p^m\rangle$. Indeed, $b_p^s$ is a power of $b_p^m$. On the other hand, $(a^kb_q^r)\cdot (b_p^s)^{\frac{n^k}{(n,s)}}\cdot (a^kb_q^r)^{-1}=a^k\cdot b_p^{mn^k}\cdot a^{-k}=b_p^m$.

\medskip

By Lemma~\ref{3.3}, $H=\langle a^kb_q^r,b^m\rangle$. We show that $H=\langle a^kb^l,b^m\rangle$ for some $l$.
If $q\geqslant 0$, then $b_q=b^{n^q}$ and we can take $l=rn^q$.
Let $q<0$. Since $(m,n)=1$, there exists an integer $t$, such that $mt\equiv r\mod (n^{-q})$. Denote $l=(r-mt)/n^{-q}$. Then, again with the help of
Lemma~\ref{3.3}, we have
$$H=\langle a^kb_q^r,b_q^m\rangle=\langle a^kb_q^{r-mt},b_q^m\rangle=\langle a^kb_q^{ln^{-q}},b_q^m\rangle=\langle a^kb^l,b^m\rangle.$$

To prove the last claim, one have to check, that $\{a^ib^j\,|\, 0\leqslant i<k, 0\leqslant j <m\}$ is
the set of representatives of the left cosets of $H$ in ${\text{\rm BS}}(1,n)$. We leave this to the reader. \hfill $\Box$

\begin{prop}\label{3.5} Let $H=\langle a^{k}b^l,b^m\rangle$ be a subgroup of ${\text{\rm BS}}(1,n)$ with $k,m>0$ and $(n,m)=1$.
Then $H$ has the presentation
$\langle x,y\,|\, x^{-1}yx=y^{n^k}\rangle$ with generators $x=a^kb^l$,  $y=b^m$.
\end{prop}

{\it Proof.} Consider the homomorphism $\psi: {\text{\rm BS}}(1,n)\rightarrow \mathbb{Z}$ introduced above.
We have $\psi(x)=k$ and $H\cap {\text{\rm ker}}\,\psi=\langle x^{-i}yx^i\,|\, i\in \mathbb{Z}\rangle $.
Thus, we have $H=\langle x^{-i}yx^i\,|\, i\in \mathbb{Z}\rangle\rtimes \langle x\rangle$.

Using the isomorphism ${\text{\rm BS}}(1,n)\cong \mathbb{Z}[\frac{1}{n}]\rtimes \mathbb{Z}$ from Lemma~\ref{3.1}, we can write
$H\cong \mathbb{Z}[\frac{m}{n^k}]\rtimes k\mathbb{Z}\cong \mathbb{Z}[\frac{1}{n^k}]\rtimes \mathbb{Z}$,
where $\mathbb{Z}$ acts on $\mathbb{Z}[\frac{1}{n^k}]$ by multiplication by $n^k$.
By Claim 3) of Lemma~\ref{3.1} we have $H\cong {\text{\rm BS}}(1,n^k)$. \hfill $\Box$


\begin{prop}\label{3.6} Let $H_1=\langle a^{k_1}b^{l_1},b^{m_1}\rangle$ and
$H_2=\langle a^{k_2}b^{l_2},b^{m_2}\rangle$ be two subgroups of ${\text{\rm BS}}(1,n)$ with $k_1,k_2,m_1,m_2>0$ and $(n,m_1)=(n,m_2)=1$. Then $H_1$ is isomorphic to $H_2$ if and only if $k_1=k_2$.
\end{prop}

{\it Proof.} If $k_1=k_2$, then $H_1\cong H_2$ by Proposition~\ref{3.5}.
This proposition also implies, that $H_i/[H_i,H_i]\cong \mathbb{Z}\times \mathbb{Z}_{n^{k_i}-1}$. So, if
$k_1\neq k_2$, then $H_1\ncong H_2$.
\hfill $\Box$

\section{The proof of the Main Theorem}

Let $\mathcal{G}$ be the subgroup of
${\text{\rm GL}}_2(\mathbb{Q})$, consisting of the matrices
$A=
\begin{pmatrix}
1 & A_{12}\\
0 & A_{22}
\end{pmatrix}
$
with $A_{12}\in \mathbb{Q}$ and $A_{22}\in \mathbb{Q}^{\ast}$.
Let $\mathcal{G}_1$ and $\mathcal{G}_2$ denote the diagonal and the unipotent subgroups of $\mathcal{G}$, i.e.
$$\mathcal{G}_1=\{A\in \mathcal{G}\,|\, A_{12}=0\},\hspace*{5mm} \mathcal{G}_2=\{A\in \mathcal{G}\,|\, A_{22}=1\}.$$
Clearly, $\mathcal{G}=\mathcal{G}_2\rtimes \mathcal{G}_1$.

For any natural $n$, let $\mathcal{H}_n$ be the subgroup of $\mathcal{G}$ consisting of the matrices $A$ with
$A_{12}\in \mathbb{Z}[\frac{1}{n}]$ and $A_{22}\in \{n^i\,|\, i\in \mathbb{Z}\}$.



\begin{lem}\label{4.1} For any natural $n\geqslant 2$, the map
$a\mapsto A=\begin{pmatrix}
1 & 0\\
0 & n
\end{pmatrix}
$,
$b\mapsto B=\begin{pmatrix}
1 & 1\\
0 & 1
\end{pmatrix}
$
can be extended to an isomorphism
$
\theta: {\text{\rm BS}}(1,n)\rightarrow \mathcal{H}_n$.
\end{lem}

{\it Proof.} The proof is easy; see Exercise 5.5 in Chapter 2 in \cite{B}. \hfill $\Box$

\medskip

We will use the following theorem of D.~Collins.

\begin{thm}\label{4.2} {\rm (\cite[Proposition~A]{C})} Let $G=\langle a,b\,|\, a^{-1}ba=b^s\rangle$ where $|s|\neq 1$. Let
$$s=\delta p_1^{e_1}p_2^{e_2}\dots p_f^{e_f},$$

where $\delta=\pm 1$ and $p_1,p_2,\dots,p_f$ are distinct primes. Then ${\text{\rm Aut}}(G)$ has presentation:
$$
\begin{array}{l}
\langle C,Q_1,Q_2,\dots, Q_f,T\,|\,\vspace*{2mm}\\ Q_i^{-1}CQ_i=C^{p_i},\,\, Q_iQ_j=Q_jQ_i,\vspace*{2mm}\\ T^2=1,\,\, TQ_i=Q_iT,\,\, T^{-1}CT=C^{-1}
\rangle,
\end{array}
$$
where $i,j=1,2,\dots ,f$.
In this presentation the automorphisms are defined by
$$
Q_i:
\begin{cases}
a \mapsto a\\
b \mapsto b^{p_i},
\end{cases}
C:
\begin{cases}
a \mapsto ab\\
b \mapsto b,
\end{cases}
T:
\begin{cases}
a \mapsto a\\
b \mapsto b^{-1}.
\end{cases}
$$
\end{thm}

\begin{prop}\label{4.3} Let $n\geqslant 2$ be a natural number.
We identify ${\text{\rm BS}}(1,n)$ with
$\mathcal{H}_n$ through the isomorphism described in Lemma~\ref{4.1}.
Let $H_1,H_2$ be two isomorphic subgroups of ${\text{\rm BS}}(1,n)$, both of finite index.
Then for every isomorphism $\varphi:H_1\rightarrow H_2$, there exists a unique matrix
$M=M(\varphi)\in \mathcal{G}$
such that $M^{-1}xM=\varphi(x)$
for every $x\in H_1$.
\end{prop}

{\it Proof.} First we prove the existence of $M(\varphi)$. By Propositions~\ref{3.4} and~\ref{3.6}, we can write $H_1=\langle a^{k}b^{l_1},b^{m_1}\rangle$ and
$H_2=\langle a^{k}b^{l_2},b^{m_2}\rangle$ for some integer $l_1,l_2$, and $k,m_1,m_2>0$, where $(n,m_1)=(n,m_2)=1$.
\!By Proposition~\ref{3.5}, $H_j$ has the presentation
$\langle x_j,y_j\,|\, x_j^{-1}y_jx_j=~y_j^{n^k}\rangle,$
where $x_j=a^{k}b^{l_j}$,  $y_j=b^{m_j}$, $j=1,2$.
After identification of elements of ${\text{\rm BS}}(1,n)$ with matrices, we have
$$x_j=\begin{pmatrix}
1 & l_j\vspace*{1mm}\\
0 & n^k
\end{pmatrix},\hspace*{2mm}
y_j=\begin{pmatrix}
1 & m_j\vspace*{1mm}\\
0 & 1
\end{pmatrix}.\eqno{(1)}
$$

Let $\varphi_0:H_1\rightarrow H_2$ be the isomorphism, such that $\varphi_0(x_1)=x_2$ and $\varphi_0(y_1)=y_2$.
Then $\varphi=\varphi_1\varphi_0$ for some $\varphi_1\in {\text{\rm Aut}}(H_1)$.
By Theorem~\ref{4.2}, ${\text{\rm Aut}}(H_1)$ is generated by the automorphisms
$$
\alpha_i:
\begin{cases}
x_1 \mapsto x_1\\
y_1 \mapsto y_1^{p_i},
\end{cases}
\beta:
\begin{cases}
x_1 \mapsto x_1y_1\\
y_1 \mapsto y_1,
\end{cases}
\gamma:
\begin{cases}
x_1 \mapsto x_1\\
y_1 \mapsto y_1^{-1},
\end{cases}
$$
$i=1,2,\dots ,f$, where $p_1,p_2,\dots,p_f$ are all prime numbers dividing $n$. Thus, it is sufficient to show
the existence of the matrices $M(\varphi_0)$, $M(\beta)$, $M(\gamma)$, and $M(\alpha_i)$, $i=1,2,\dots,f$.

First we prove the existence of $M(\varphi_0)$. We shall find
$M(\varphi_0)\in
\mathcal{G}$,
such that
$$
\begin{array}{l}
x_1\cdot M(\varphi_0)=M(\varphi_0)\cdot \varphi_0(x_1),\vspace*{2mm}\\
y_1\cdot M(\varphi_0)=M(\varphi_0)\cdot \varphi_0(y_1).
\end{array}
$$
Using (1), one can compute that
$$M(\varphi_0)=
\begin{pmatrix}
1 & \frac{l_1m_2-l_2m_1}{m_1(n^k-1)}\vspace*{2mm}\\
0 & \frac{m_2}{m_1}
\end{pmatrix}.\eqno{(2)}
$$
Similarly, we get
$$
M(\alpha_i)=
\begin{pmatrix}
1 & \frac{l_1(p_i-1)}{n^k-1}\vspace*{2mm}\\
0 & p_i
\end{pmatrix},\hspace*{2mm}
M(\beta)=
\begin{pmatrix}
1 & \frac{-m_1}{n^k-1}\vspace*{2mm}\\
0 & 1
\end{pmatrix},\hspace*{2mm}
M(\gamma)=
\begin{pmatrix}
1 & \frac{-2l_1}{n^k-1}\vspace*{2mm}\\
0 & -1
\end{pmatrix}.\eqno{(3)}
$$
The uniqueness of $M$ follows from the triviality of the centralizer of $H_1$ in $\mathcal{G}$; the later is easy to check.  \hfill $\Box$


\begin{lem}\label{4.4} {\rm 1)} Let $\varphi:H\rightarrow H'$ be an isomorphism between subgroups of finite index in ${\text{\rm BS}}(1,n)$ and let $K$ be a subgroup of finite index in $H$. Then $M(\varphi|_K)=M(\varphi)$.

{\rm 2)} Let $\varphi_1: H_1\rightarrow H_1'$ and
$\varphi_2: H_2\rightarrow H_2'$ be two isomorphisms between subgroups of finite index in ${\text{\rm BS}}(1,n)$.
Suppose that $[\varphi_1]=[\varphi_2]$ in
${\text{\rm Comm(BS}}(1,n))$. Then $M(\varphi_1)=M(\varphi_2)$.
\end{lem}

{\it Proof.} 1) For every $x\in K$ we have
$M(\varphi|_K)^{-1}xM(\varphi|_K)=\varphi|_K(x)=\varphi(x)=M(\varphi)^{-1}xM(\varphi)$ and the claim
follows from the uniqueness of $M$.

2) By Lemmas~\ref{2.4} and~\ref{2.5}, we have $\varphi_1|_{H_1\cap H_2}=\varphi_2|_{H_1\cap H_2}$.
Claim 1) implies that $M(\varphi_1)=M(\varphi_1|_{H_1\cap H_2})=M(\varphi_2|_{H_1\cap H_2})=M(\varphi_2)$. \hfill $\Box$

\medskip

This enables to define $M$ on commensurator classes: $M([\varphi]):=M(\varphi)$.

\begin{thm}\label{4.5} For every natural $n\geqslant 2$, the map $\Psi:{\text{\rm Comm(BS}}(1,n))\rightarrow
\mathcal{G}$
given by
$[\varphi]\mapsto M([\varphi])$ is an isomorphism.
\end{thm}

{\it Proof.} 1) First we prove that $\Psi$ is a homomorphism.
Let $\varphi_1: H_1\rightarrow H_2$, $\varphi_2: H_3\rightarrow H_4$ be two isomorphisms between subgroups of finite index in ${\text{\rm BS}}(1,n)$.
We shall show that $M([\varphi_1])M([\varphi_2])=M([\varphi_1\varphi_2])$. Write $\varphi_1\varphi_2=\sigma\tau$, where $\sigma$ is the restriction of $\varphi_1$ to $\varphi_1^{-1}(H_2\cap H_3)$ and $\tau$ is the
restriction of $\varphi_2$ to $H_2\cap H_3$:
$$\varphi_1^{-1}(H_2\cap H_3)\overset{\sigma}{\longrightarrow} (H_2\cap H_3)\overset{\tau}{\longrightarrow} \varphi_2(H_2\cap H_3).$$

For $x\in \varphi_1^{-1}(H_2\cap H_3)$ we have $(\varphi_1\varphi_2)(x)=\tau((\sigma(x)))=M(\tau)^{-1}M(\sigma)^{-1}x
M(\sigma)M(\tau)$. Hence, $M(\varphi_1\varphi_2)=M(\sigma)M(\tau)=M(\varphi_1)M(\varphi_2)$ and the claim follows.

2) The injectivity of $\Psi$ trivially follows from the definition of $M([\varphi])$.

3) Now we prove that $\Psi$ is a surjection. By specializing parameters in (2) and (3), we will obtain some matrices in ${\text{\rm im}}\Psi$.
Taking $l_1=m_2$ and $l_2=m_1$ in $M(\varphi_0)$, we get the matrix $$D(\begin{matrix}
\frac{m_1}{m_2}\end{matrix})=\begin{pmatrix}
1 & 0\\
0 & \frac{m_1}{m_2}
\end{pmatrix}$$ with $m_1,m_2>0$, $(m_1,n)=(m_2,n)=1$.
Taking $l_1=0$ in $M(\alpha_i)$ and in $M(\gamma)$, and taking $m_1=1$ in $M(\beta)$, we get the matrices
$$D(p_i)=\begin{pmatrix}
1 & 0\\
0 & p_i,
\end{pmatrix},\hspace*{4mm}
D(-1)=\begin{pmatrix}
1 & 0\\
0 & -1
\end{pmatrix},\hspace*{4mm}
T(k)=\begin{pmatrix}
1 & \frac{1}{n^k-1}\\
0 & 1
\end{pmatrix},\hspace*{2mm}
k>0.
$$
The matrices $D(\frac{m_1}{m_2})$, $D(p_i)$ and $D(-1)$ generate the subgroup $\mathcal{G}_1$
in the image of $\Psi$.

So, it suffices to show that $\mathcal{G}_2$
is contained in ${\text{\rm im}}\Psi$.
Since the additive group of $\mathbb{Q}$ is generated by $\mathbb{Z}[\frac{1}{n}]$ and all numbers $\frac{1}{s}$ with $(s,n)=1$,
it suffices to show that the matrices
$
\begin{pmatrix}
1 & q\\
0 & 1
\end{pmatrix}
$
with $q\in \mathbb{Z}[\frac{1}{n}]$
and the matrices
$\begin{pmatrix}
1 & \frac{1}{s}\\
0 & 1
\end{pmatrix}
$
with $(s,n)$=1 are contained in the image of~$\Psi$.
The first follows from the fact that
the  group of the commensurator classes of inner
automorphisms of ${\text{\rm BS}}(1,n)$ is mapped, under $\Psi$, onto $\mathcal{H}_n$.
The second follows from the formula
$\begin{pmatrix}
1 & \frac{1}{s}\\
0 & 1
\end{pmatrix}=(T(\phi(s)))^t,$
where $\phi$ is the Euler function and $t$ is the natural number such that $n^{\phi(s)}-1=st$. \hfill $\Box$

\section{Appendix: Commensurators and quasi-isometries}

Let $X$ and $Y$ be two metric spaces. A map $f: X\rightarrow Y$ is called a (coarse) {\it quasi-isometry} between $X$ and $Y$ if there are some constants $K,C,C_0 > 0$, such that the following holds:

\bigskip

1. $K^{-1} d_X(x_1,x_2)-C\leqslant d_Y(f(x_1),f(x_2))\leqslant K d_X(x_1,x_2)+C$ for all $x_1,x_2\in X$.

2. The $C_0$-neighborhood of $f(X)$ coincides with $Y$.

\medskip

There is always a coarse inverse of $f$, a quasi-isometry $g :Y\rightarrow X$ such
that $f\circ g$ and $g\circ f$ are a bounded distance from the identity maps in the sup
norm; these bounds, and the quasi-isometry constants for $g$, depend only on
the quasi-isometry constants of $f$.

\begin{defn}\label{5.1} {\rm
Let $X$ be a metric space. Two quasi-isometries $f$ and $g$ from $X$ to itself are
considered equivalent if there exists a number $M >0$ such that $d(f(x),g(x))\leqslant M$ for all
$x\in X$.
Let ${\text{\rm QI}}(X)$ be the set of equivalence
classes of quasi-isometries from $X$ to itself. Composition of quasi-isometries
gives a well-defined group structure on ${\text{\rm QI}}(X)$. The group ${\text{\rm QI}}(X)$
is called the {\it quasi-isometry group} of $X$.}

\end{defn}

Let $G$ be a group with a finite generating set $S$.
For $g\in G$ denote by $|g|$ the minimal $k$, such that $g=s_1s_2\dots s_k$, where $s_1,s_2,\dots,s_k\in S\cup S^{-1}$.
We consider $G$ as a metric space with the word metric with respect to $S$: $d(x,y)=|x^{-1}y|$ for $x,y\in G$.
For a finitely generated group $G$, the group ${\text{\rm QI}}(G)$ is well defined and does not depend
on a choice of a finite generating set $S$.

It is well known that there is a natural homomorphism $\Lambda: {\text{\rm Comm}}(G)\rightarrow {\text{\rm QI}}(G)$.
This homomorphism is defined by the following rule. Let $\varphi:H\rightarrow H'$ be
an isomorphism between two finite index subgroups of $G$. We choose a right transversal $T$ for $H$ in $G$ with $1\in T$.
First we define a map $f_{\varphi}:G\rightarrow G$ by the rule $f_{\varphi}(ht):=\varphi(h)$ for every $h\in H$ and $t\in T$. Clearly, $f_{\varphi}$ is a quasi-isometry. Then we set $\Lambda([\varphi]):=[f_{\varphi}]$.

\begin{prop}\label{5.2}  Let $G$ be a finitely generated group with the unique root property. Then $\Lambda:{\text{\rm Comm}}(G)\rightarrow {\text{\rm QI}}(G)$ is an embedding.
\end{prop}

{\it Proof.} We will use notation introduced before this lemma. Suppose that $[f_{\varphi}]=[{\text{\rm id}}_{|G}]$.
Then there is a constant $M>0$, such that $d(f_{\varphi}(x),x)\leqslant M$ for every $x\in G$.
Let $h\in H$. Then for every integer $n$ holds: $|h^{-n}\varphi(h^n)|=d(\varphi(h^n),h^n)\leqslant M$.
Since $G$ is finitely generated, the $M$-ball in $G$ centered at 1 is finite.
Hence, there exist distinct $n,m$ such that $h^{-n}\varphi(h^n)=h^{-m}\varphi(h^m)$.
Then $h^{n-m}=(\varphi(h))^{n-m}$ and so $h=\varphi(h)$ by the unique root property.
Hence $[\varphi]=1$ and the injectivity of $\Lambda$ is proved. \hfill $\Box$

\begin{cor}\label{5.3} The group ${\text{\rm Comm}}({\text{\rm BS}}(1,n))$ naturally embeds in ${\text{\rm QI}}({\text{\rm BS}}(1,n))$.
\end{cor}

\section{Acknowledgements}
The author thanks the MPIM at Bonn for its support and excellent working conditions during the fall 2010, while this research was finished. He is gratefully thankful to Nadine Hansen for hearing the first version of this paper.


\begin{thebibliography}{AA}

\bibitem{BB} L. Bartholdi, O. Bogopolski, {\it On abstract commensurators of groups},
J. Group Theory, {\bf 13}, (6) (2010), 903-922.

\bibitem{BS} G. Baumslag G, D. Solitar, {\it Some two-generator, one-relator non-Hopfian
groups}, Bull. Amer. Math. Soc., {\bf 68} (1962), 199-201.

\bibitem{B} O. Bogopolski, {\it Introduction to group theory}, EMS: Z\"{u}rich, 2008.

\bibitem{C} D. Collins, {\it The automorphism towers of some one-relator groups}, Proc. London Math. Soc.,
{\bf 36}, (3) (1978), 480-493.

\bibitem{FM} B. Farb and L. Mosher (appendix by D. Cooper), {\it A rigidity theorem for the solvable
Baumslag-Solitar groups}, Inventiones, {\bf 131},  (2) (1998), 419-451.

\bibitem{W} K. Whyte, {\it The large scale geometry of the higher Baumslag -- Solitar groups}, GAFA, {\bf 11} (2001),
1327-1343.



\end{thebibliography}
\end{document}